\documentclass[12pt,twoside]{article}

\usepackage{amsfonts,amsmath,amssymb,amsthm,amscd,graphicx,color,accents}

\usepackage[colorlinks=true,linkcolor=blue,citecolor=blue,urlcolor=blue]{hyperref}

\numberwithin{equation}{section}

\theoremstyle{plain}
\newtheorem{theorem}{Theorem}[section]
\newtheorem{proposition}[theorem]{Proposition}
\newtheorem{lemma}[theorem]{Lemma}
\newtheorem{corollary}[theorem]{Corollary}
\newtheorem{remark}[theorem]{Remark}

\newcommand{\Pro}{\noindent\textit{Proof.}\ \ }
\def \dd {\,\textrm{d}}

\def\E{\mathbb{E}\hskip 1pt}
\def\Pr {\mathbb{P}}
\def\R{\mathbb{R}}
\def\N{\mathbb{N}}
\def\tr{{\rm{tr}}\hskip 1pt}
\def\diag{{\rm{diag}}}
\def\bX{{\boldsymbol{X}}}
\def\bk{{\boldsymbol{k}}}
\def\bl{{\boldsymbol{l}}}
\def\bm{{\boldsymbol{m}}}
\def\bn{{\boldsymbol{n}}}

\begin{document}

\title{\bf\Large Product Inequalities for Multivariate Gaussian, Gamma, and Positively Upper Orthant Dependent Distributions}

\author{{Dominic Edelmann}\thanks{
Division of Biostatistics, German Cancer Research Center, Im Neuenheimer Feld 280, 69120 Heidelberg, Germany.  E-mail: \href{mailto:dominic.edelmann@dkfz-heidelberg.de}{dominic.edelmann@dkfz-heidelberg.de}.
\endgraf
\ $^\mathparagraph$Corresponding Author, E-mail: \href{mailto:dominic.edelmann@dkfz-heidelberg.de}{dominic.edelmann@dkfz-heidelberg.de}. 
\endgraf 
\ $^\dag$Department of Statistics, Pennsylvania State University, University Park, PA 16802, U.S.A. E-mail: \href{mailto:richards@stat.psu.edu}{richards@stat.psu.edu}.
\endgraf
\ $^\ddag$University of Applied Sciences, D-55411 Bingen, Germany.  E-mail: \href{mailto:thomas.royen@t-online.de}{thomas.royen@t-online.de}.
\endgraf
\endgraf
\endgraf
}\, $^{,\mathparagraph}$, \ {Donald Richards}$^\dag$,
 \ and {Thomas Royen}$^\ddag$
 \endgraf
}

\maketitle

\begin{abstract}
The Gaussian product inequality is an important conjecture concerning the moments of Gaussian random vectors.  While all attempts to prove the Gaussian product inequality in full generality have been unsuccessful to date, numerous partial results have been derived in recent decades and we provide here further results on the problem. Most importantly, we establish a strong version of the Gaussian product inequality for multivariate gamma distributions in the case of nonnegative correlations, thereby extending a result recently derived by Genest and Ouimet \cite{Ouimet}. Further, we show that the Gaussian product inequality holds with nonnegative exponents for all random vectors with positive components whenever the underlying vector is positively upper orthant dependent. Finally, we show that the Gaussian product inequality with negative exponents follows directly from the Gaussian correlation inequality.
\end{abstract}

\section{Introduction}
\label{sec:introduction}

Let $d$ be a positive integer, and let $\bX = (X_1,\ldots,X_d)$ be a random vector which has a multivariate Gaussian distribution with probability density function,
$$
(2\pi)^{-n/2} |\Sigma|^{-1/2} \exp(-\tfrac12 x'\Sigma^{-1} x), \qquad x \in \R^d,
$$
with a nonsingular covariance matrix $\Sigma$.  We refer to the random vector $\bX$ as having a \textit{centered} Gaussian distribution because $\E(\bX) = 0$, and we write $\bX \sim \mathcal{N}_d(0,\Sigma)$.  

The Gaussian product inequality (GPI) conjecture states that, for any centered Gaussian random vector $\bX = (X_1,\ldots,X_d)$ and any $n \in \mathbb{N}$, there holds the inequality, 
\begin{equation} 
\label{eq:classgpi}
\E \Big(\prod_{j=1}^d X_j^{2n}\Big) \geq \prod_{j=1}^d \E(X_j^{2n}). 
\end{equation}
We refer readers to Frenkel \cite{Frenkel} and Genest and Ouimet \cite{Ouimet} for details of the history, motivation, and literature on this inequality.  

Several generalizations of \eqref{eq:classgpi} have been studied recently.  Li and Wei \cite{Liwei} considered, as an extension of \eqref{eq:classgpi}, conditions such that 
\begin{equation}
\label{weak_gpi}
\E \Big(\prod_{j=1}^d |X_j|^{n_j}\Big) \ge \prod_{j=1}^d \E\big(|X_j|^{n_j}\big),
\end{equation}
for arbitrary $n_1,\ldots,n_d > 0$.  Wei \cite{Wei} derived hypotheses under which the inequality
\begin{equation}
\label{strong_gpi}
\E \Big(\prod_{j=1}^d |X_j|^{n_j}\Big) \ge \E\Big(\prod_{j \in I} |X_j|^{n_j}\Big) \cdot \E\Big(\prod_{j \in I^c} |X_j|^{n_j}\Big),
\end{equation}
holds for all index sets $I \subset \{1,\ldots,d\}$, where $I^c = \{1,\ldots,d\} \setminus I$. To distinguish between these inequalities, we call \eqref{weak_gpi} the {\it weak form} and \eqref{strong_gpi} the {\it strong form} of the GPI.

Russell and Sun \cite{Russell} recently related the GPI to a class of combinatorial inequalities, and thereby established numerous cases of the GPI for $d=3$.  One of the results obtained by Russell and Sun \cite{Russell} is derived here by different methods, and we present that result in Corollary \ref{cor_normal_strong_GPI}.  The approach by way of combinatorial inequalities is noteworthy because it is also shown in \cite{Russell} to lead to new inequalities for the bivariate Gaussian distributions.  

The weak form of the GPI was established by Frenkel \cite{Frenkel} for $n_1 = \cdots = n_d = 2$ and arbitrary $d$, and by Lan, \textit{et al.}~\cite{Lan} for $d = 3$ and integer exponents $n_1$, $n_2$, and $n_3$ with equality between at least two exponents.  Genest and Ouimet \cite{Ouimet} developed recently a novel and far-reaching approach to the GPI, proving \eqref{weak_gpi} for arbitrary $d$ with nonnegative even integer exponents $n_j$ when the covariance matrix $\Sigma$ is \textit{completely positive}, i.e., $\Sigma = CC'$ where $C = (c_{ij})$ is a $d \times d$ matrix with $c_{ij} \ge 0$ for all $i,j=1,\ldots,d$.  Nevertheless, it is still unknown whether the weak form of the GPI \eqref{weak_gpi} is valid for general $\Sigma$.  

On the other hand, the strong form of the GPI \eqref{strong_gpi} fails even for $d=3$ with $n_1 = n_2 = n_3 = 2$. Consider a Gaussian random vector, $\bX \sim \mathcal{N}_3(0,\Sigma)$, with $\Sigma = (\sigma_{ij})$.  By the Isserlis-Wick formula \cite{Kan} or by using moment-generating functions, we obtain 
$$
\E(X_1^2 X_2^2 X_3^2) - \E(X_1^2 X_2^2) \E(X_3^2) = 2 \big(\sigma_{11} \sigma_{23}^2 + 4 \sigma_{12} \sigma_{13} \sigma_{23} + \sigma_{13}^2 \sigma_{22}\big).
$$
If we set $\sigma_{ii} = 1$, $i=1,2,3$, and $\sigma_{13} = \sigma_{23} = \rho$ with $0 < |\rho| < 1/2$ then with $\sigma_{12} \in (2 \rho^2 -1,- 1/2)$, the matrix $\Sigma$ is positive definite and yet 
$$
\E(X_1^2 X_2^2 X_3^2) - \E(X_1^2 X_2^2) \E(X_3^2) = 4 \rho^2 (1 + 2 \sigma_{12}) < 0;
$$
concrete examples are $(\sigma_{12},\rho) = (-0.6,0.4)$ and $(\sigma_{12},\rho) = (-0.75,0.3)$.  

Wei \cite{Wei} showed, however, that the strong form of the GPI holds for the case in which all exponents $n_1,\ldots,n_d$ are negative.  It is also obvious that the strong form \eqref{strong_gpi} holds if $|\bX| := (|X_1|,\ldots,|X_d|)$, the vector of absolute values, is \textit{associated}, i.e., if 
$\mbox{Cov} \big(f(|\bX|),g(|\bX|)\big) \geq 0$ for all component-wise non-decreasing functions $f,g : \R_+^n \to \R$ \cite{esary}.  Thus, for a centered Gaussian random vector $\bX$, if $|\bX|$ is associated then it follows immediately that the strong form of the GPI holds.  In particular, if the vector $|\bX|$ is \textit{multivariate totally positive of order} $2$, denoted $\mbox{MTP}_2$ (cf., \cite{Karlin1}) then, as the $\mbox{MTP}_2$ property implies associatedness, it follows that the strong form of the GPI holds.  

Moreover, for Gaussian vectors $\bX$, the $\mbox{MTP}_2$ property of its absolute values $|\bX|$ can be characterized explicitly in terms of the covariance matrix $\Sigma$. For this purpose (and in the sequel) we call a diagonal matrix $S = \mbox{diag}(s_1,\ldots,s_d)$ a \textit{sign matrix} if $s_j = \pm 1$ for all $j =1,\ldots,d$.  It was proved by Karlin and Rinott \cite{Karlin2} that for $X \sim \mathcal{N}_d(0,\Sigma)$ the vector of absolute values, $|\bX|$, is $\mbox{MTP}_2$ if and only if there exists a sign matrix $S$ such that all off-diagonal entries of $- S \Sigma^{-1} S$ are nonnegative; hence the strong form of the GPI holds for that class of covariance matrices.

In this article, we derive new and more general hypotheses under which the weak form of the GPI \eqref{weak_gpi} and the strong form of the GPI \eqref{strong_gpi} hold.  We extend the results of Genest and Ouimet \cite{Ouimet} in several directions, one of which is a proof of the strong form of the GPI \eqref{strong_gpi} for nonnegative correlations, i.e., for any covariance matrix $\Sigma = (\sigma_{ij})$ with $\sigma_{ij} \ge 0$ for all $1 \le i < j \le d$.  Additionally, we show that the weak form of the GPI and the strong form of the GPI follow from the properties of positive upper orthant dependence (PUOD) and strongly positive upper orthant dependence (SPUOD), respectively.  Finally, we apply the Gaussian correlation inequality (Royen \cite{Royen}) to obtain in Section \ref{sec_negative_exponents} an alternative and succinct proof of the strong form of the GPI for negative exponents, derived originally by Wei \cite{Wei}; further, we show that this result extends to the multivariate gamma distributions.

\section{The strong form of the GPI for nonnegative correlations}
\label{sec_corr_nonneg}

Genest and Ouimet \cite{Ouimet} established the weak form of the GPI \eqref{weak_gpi} for the multivariate normal distribution $\mathcal{N}_d(0,\Sigma)$ and with even integers $n_1,\ldots,n_d$ and completely positive covariance matrix $\Sigma$, i.e., $\Sigma = CC'$, where $C = (c_{ij})$ is a matrix with nonnegative entries $c_{ij} \geq 0$.  In Theorem \ref{thm_nonneg_corr}, we extend this result in three directions. First and most importantly, we extend the result in \cite{Ouimet} to the case of nonnegative correlations, where $\Sigma = (\sigma_{ij})$ is such that $\sigma_{ij} \ge 0$ for all $1 \le i < j \le d$.  For $d \geq 5$, the assumption of nonnegative correlations is known to be less restrictive than complete positivity \cite{Gray}. 

On the other hand, a famous counterexample of \v{S}id\'ak \cite{Sidak} established the existence of Gaussian random vectors $\bX = (X_1,X_2,X_3)$ with completely positive covariance matrices and for which the vector $|\bX| = (|X_1|,|X_2|,|X_3|)$ of absolute values is not positively upper orthant dependent (PUOD) and hence not associated.  Hence the result of Genest and Ouimet and the more general result presented here both extend the strong form of the GPI \eqref{strong_gpi} beyond the straightforward case in which the vector of absolute values is associated.  

The second direction in which Theorem \ref{thm_nonneg_corr} extends results known hitherto is that the we obtain the strong form \eqref{strong_gpi}, hence also the weak form \eqref{weak_gpi}.

Third, in considering the case of even exponents, the weak form of the GPI \eqref{weak_gpi} and the strong form \eqref{strong_gpi} each correspond to inequalities for special cases of the multivariate gamma distributions.  Precisely, the $d$-dimensional gamma distribution (in the sense of Krishnamoorthy and Parthasarathy \cite{Krishna}) may be defined by means of its moment-generating function.  Denote by $I_d$ the identity matrix of order $d$ and, for sufficiently small $t_1,\ldots,t_d \in \mathbb{R}$, define $T = \diag(t_1,\ldots,t_d)$.  Then we say that $\bX = (X_1,\ldots,X_d)$ has a multivariate gamma distribution with a not necessarily integer ``degree-of-freedom parameter'' $2 \alpha$ and positive semidefinite matrix parameter $\Sigma$, written $\bX \sim {\hbox{Gamma}\hskip1pt}_d(\alpha,\Sigma)$, if the moment-generating function of $\bX$ is 
\begin{equation} 
\label{eq:lpgamma}
\E \exp\Big(\sum_{j=1}^d t_j X_j\Big) = \det(I_d - \Sigma T)^{-\alpha}.
\end{equation}
This $d$-dimensional gamma distribution is also known as the \textit{Wishart-Gamma distribution} since it was derived originally as the distribution of one-half of the diagonal entries of a $W_d(2 \alpha,\Sigma)$-Wishart distributed random matrix with $2 \alpha \in \mathbb{N} \cup (d-1,\infty)$. In this regard, it is remarkable that the ${\hbox{Gamma}\hskip1pt}_d(\alpha,\Sigma)$ distribution also exists for all values $2 \alpha > [(d-1)/2]$, where $[u]$ denotes the integer part of $u \in \R$; see \cite{Royen_gamma}.  

The ${\hbox{Gamma}\hskip1pt}_d(\frac{1}{2},\Sigma)$ distribution is infinitely divisible (i.e., the ${\hbox{Gamma}\hskip1pt}_d(\alpha,\Sigma)$ distribution exists for all $\alpha > 0$) if and only if the ${\hbox{Gamma}\hskip1pt}_d(\frac{1}{2},\Sigma)$ distribution is multivariate totally positive of order $2$ ($\mbox{MTP}_2$); see Bapat \cite{Bapat}.  For example, one can show that ${\hbox{Gamma}\hskip1pt}_d(\frac{1}{2},\Sigma)$ is infinitely divisible if $\Sigma = (\sigma_{ij})$ is of ``structure $\ell$'' \cite{Royen_onefac}, i.e., $\sigma_{ij} = a_i a_j (\sigma_{ii} \sigma_{jj})^{1/2}$ with $|a_j| < 1$ for all $j = 1,\ldots,d$, or if $\Sigma^{-1}$ is of ``tree-type''; see \cite{Royen_tree} for details. 

If $(X_1,\ldots,X_d) \sim N_d(0,\Sigma)$ then $\tfrac{1}{2} (X_1^2,\ldots,X_d^2) \sim {\hbox{Gamma}\hskip1pt}_{d}(\frac{1}{2},\Sigma)$ distribution.  Consequently, in the case of even exponents, as considered in \cite{Ouimet}, the weak form \eqref{weak_gpi} and the strong form \eqref{strong_gpi} of the GPI intrinsically are inequalities on the ${\hbox{Gamma}\hskip1pt}_{d}(\frac{1}{2},\Sigma)$ distribution.  Therefore it is natural to extend these inequalities to the more general multivariate gamma distributions having moment-generating function \eqref{eq:lpgamma}.


\begin{theorem}
	\label{thm_nonneg_corr}
	Let $(X_1,\ldots,X_d) \sim {\hbox{Gamma}\hskip1pt}_d(\alpha,\Sigma)$, where $\Sigma = (\sigma_{ij})$ is positive semidefinite.  Suppose there exists a sign matrix $S$ such that all the elements in $S \Sigma S$ are nonnegative, i.e., $\sigma_{ij} s_i s_j \ge 0$ for all $1 \le i < j \le d$.  Then for all subsets $I \subset \{1,\ldots,d\}$, and for all nonnegative integers $n_1,\ldots,n_d$, there holds the strong form of the GPI,
	\begin{equation}
	\label{strong_gpi_gamma}
	\E \Big(\prod_{j=1}^d X_j^{n_j}\Big) \ge \E\Big(\prod_{j \in I} X_j^{n_j}\Big) \cdot \E\Big(\prod_{j \in I^c} X_j^{n_j}\Big).
	\end{equation}
\end{theorem}

\Pro
Since the moment-generating function \eqref{eq:lpgamma} is invariant under the transformation $\Sigma \to S \Sigma S$ we can, without loss of generality, assume that $\sigma_{ij} \ge 0$ for all $1 \le i < j \le d$.  Moreover, by permuting the coordinates of $\bX$, we may also assume that $I=\{1,\ldots,p\}$ where $1 \le p \le d-1$.

With sufficiently small $t_1,\ldots,t_d \in [0,1)$ and $T = \diag(t_1,\ldots,t_d)$, the moment-generating function of the ${\hbox{Gamma}\hskip1pt}_d(\alpha,\Sigma)$ distribution is
\begin{equation}
\label{eq_mgf_full_X}
\E \exp\Big(\sum_{j=1}^d t_j X_j\Big) = \det(I_d - \Sigma T)^{-\alpha} = \det(I_d - T^{1/2} \Sigma T^{1/2})^{-\alpha}.
\end{equation}
Denote by $\varepsilon_1,\ldots,\varepsilon_d$ the eigenvalues of the matrix $T^{1/2} \Sigma T^{1/2}$.  For sufficiently small $t_1,\ldots,t_d \in [0,1)$, we have $\varepsilon_1,\ldots,\varepsilon_d \in [0,1)$.  Then we have 
\begin{align*}
\det(I_d - T^{1/2} \Sigma T^{1/2})^{-\alpha} &= \exp\big(-\alpha \log \det(I_d - T^{1/2} \Sigma T^{1/2}\big) \\
&= \exp\Big(-\alpha \sum_{j=1}^d \log (1 - \varepsilon_j)\Big).
\end{align*}
Inserting into this sum the series expansions,  
$$
-\log (1 - \varepsilon_j) = \sum_{n=1}^\infty \frac{\varepsilon_j^n}{n},
$$
$j=1,\ldots,n$, and interchanging the order of summation, we obtain 
\begin{align}
\label{eq_mgf_X}
\det(I_d - \Sigma T)^{-\alpha} &= \exp\Big(\alpha \sum_{j=1}^d \sum_{n=1}^\infty \frac{\varepsilon_j^n}{n}\Big) \nonumber \\
&= \exp\Big(\alpha \sum_{n=1}^\infty \frac{1}{n} \sum_{j=1}^d \varepsilon_j^n\Big) = \exp\Big(\alpha \sum_{n=1}^\infty \frac{1}{n} \tr[(\Sigma T)^n]\Big).
\end{align}

Next, we partition $\Sigma$ and $T$ into block matrices, 
$$
\Sigma = 
\begin{pmatrix}
\Sigma_{11} & \Sigma_{12} \\ \Sigma_{21} & \Sigma_{22}
\end{pmatrix}, \qquad 
T = \begin{pmatrix}
T_1 & 0 \\ 0 & T_2
\end{pmatrix},
$$
where $\Sigma_{11}$ and $T_1 = \diag(t_1,\ldots,t_p)$ are $p \times p$, $\Sigma_{12} = \Sigma_{21}'$ is $p \times (d-p)$, and $\Sigma_{22}$ and $T_2 = \diag(t_{p+1},\ldots,t_d)$ are $(d-p) \times (d-p)$.  Then, 
\begin{equation}
\label{eq_TSigma}
T^{1/2} \Sigma T^{1/2} = \begin{pmatrix}
T_1^{1/2} \Sigma_{11} T_1^{1/2} & T_1^{1/2} \Sigma_{12} T_2^{1/2} \\ T_2^{1/2} \Sigma_{21} T_1^{1/2} & T_2^{1/2} \Sigma_{22} T_2^{1/2}
\end{pmatrix}
\end{equation}
Let
$$
A = \begin{pmatrix}
A_{11} & A_{12} \\ A_{21} & A_{22}
\end{pmatrix}
$$
be a symmetric matrix which has been partitioned similarly to $\Sigma$.  By induction on $n$, we find that 
\begin{equation}
\label{eq_matrix_power}
A^n = \begin{pmatrix}
A_{11}^n + Q_{11,n}(A) & Q_{12,n}(A) \\
Q_{21,n}(A) & A_{22}^n + Q_{22,n}(A)
\end{pmatrix}
\end{equation}
where each matrix $Q_{ij,n}(A)$ is a homogeneous polynomial in $A$ with nonnegative coefficients; for instance, 
\begin{align*}
Q_{11,2}(A) = A_{12} A_{21}, \ Q_{12,2}(A) = [Q_{21,2}(A)]' = A_{11} A_{12} + A_{12} A_{22}, \ 
Q_{22,2} = A_{21} A_{12}.
\end{align*}
Denote by $\mathbb{N}_0$ the set of nonnegative integers.  Applying \eqref{eq_matrix_power} to \eqref{eq_TSigma}, and taking traces, we obtain 
\begin{align}
\label{eq_tr_SigmaT_n}
\frac{1}{n}\tr[(\Sigma T)^n] &= \frac{1}{n} \tr\big[(T^{1/2} \Sigma T^{1/2})^n\big] \nonumber \\
&= \frac{1}{n}\tr\big[(T_1^{1/2} \Sigma_{11} T_1^{1/2})^n\big] + \frac{1}{n} \tr\big[(T_2^{1/2} \Sigma_{22} T_2^{1/2})^n\big] + \frac{1}{n}\sum_{i=1}^2 Q_{ii,n}(T^{1/2} \Sigma T^{1/2}) \nonumber \\
&= \frac{1}{n}\tr[(\Sigma_{11} T_1)^n] + \frac{1}{n}\tr[(\Sigma_{22} T_2)^n] + \sum_{\substack{\bn \in \mathbb{N}_0^d, \\ n_1+\cdots+n_d=n}} c_\bn t_1^{n_1} \cdots t_d^{n_d},
\end{align}
where each $c_\bn$ is a polynomial in the entries of $\Sigma$ and $c_{(0,\ldots,0)} = 0$.  It is evident from \eqref{eq_matrix_power} that the coefficients of each $c_\bn$ are nonnegative; therefore, since the entries of $\Sigma$ are nonnegative, we obtain $c_\bn \ge 0$ for all $\bn$.  

Substituting \eqref{eq_tr_SigmaT_n} into \eqref{eq_mgf_X}, we obtain 
\begin{align}
&\det(I_d - \Sigma T)^{-\alpha} \nonumber \\
&= \exp\Big(\alpha \sum_{n=1}^\infty \frac{1}{n} \tr[(\Sigma_{11} T_1)^n]\Big) \exp\Big(\alpha \sum_{n=1}^\infty \frac{1}{n} \tr[(\Sigma_{22} T_2)^n]\Big) \exp\Big(\sum_{\bn \in \mathbb{N}_0^d} \alpha c_\bn t_1^{n_1} \cdots t_d^{n_d}\Big) \nonumber \\
&= \det(I_p - \Sigma_{11} T_1)^{-\alpha} \det(I_{d-p} - \Sigma_{22} T_2)^{-\alpha} \cdot \exp\Big(\sum_{\bn \in \mathbb{N}_0^d} \alpha c_\bn t_1^{n_1} \cdots t_d^{n_d}\Big). \label{eq:det}
\end{align}

Next, the Maclaurin expansion of the exponential function leads to 
\begin{equation}
\label{eq_exp3}
\exp\Big(\sum_{\bn \in \mathbb{N}_0^d } \alpha c_\bn t_1^{n_1} \cdots t_d^{n_d}\Big) = \sum_{\bm \in \mathbb{N}_0^d} b_\bm t_1^{m_1} \cdots t_d^{m_d},
\end{equation}
Since $c_{(0,\ldots,0)} = 0$ and $c_\bn \ge 0$ for all $\bn$ then $b_{(0,\ldots,0)} = 1$ and $b_\bm \ge 0$ for all $\bm$.  

Applying \eqref{eq_mgf_full_X} and \eqref{eq_exp3} to \eqref{eq:det}, we obtain 
\begin{align*}
\sum_{\bn \in \mathbb{N}_0^d} & \Big(\E\prod_{j=1}^d \frac{X_j^{n_j}}{n_j!}\Big) t_1^{n_1} \cdots t_d^{n_d} \\
&= \Big(\sum_{\bk \in \mathbb{N}_0^p} \E\prod_{j=1}^p \frac{(t_j X_j)^{k_j}}{k_j!}\Big) \cdot \Big(\sum_{\bl \in \mathbb{N}_0^{d-p}} \E\prod_{j=p+1}^d \frac{(t_j X_j)^{l_j}}{l_j!}\Big) \cdot \Big(\sum_{\bm \in \N_0^d} b_\bm t_1^{m_1} \cdots t_d^{m_d}\Big).
\end{align*}
Collecting terms in $\bn$ on the right-hand side of the above expression, we obtain 
$$
\sum_{\bn \in \mathbb{N}_0^d} \Big(\E\prod_{j=1}^d  \frac{X_j^{n_j}}{n_j!}\Big) t_1^{n_1} \cdots t_d^{n_d} = \sum_{\bn \in \mathbb{N}_0^d} \bigg[\sum_{(\bk,\bl)+\bm=\bn} b_\bm \Big(\E\prod_{j=1}^p \frac{X_j^{k_j}}{k_j!}\Big) \cdot \Big(\E\prod_{j=p+1}^d \frac{X_j^{l_j}}{l_j!}\Big)\bigg] \prod_{j=1}^d t_j^{n_j}.
$$
Next, we decompose the inner summation into terms corresponding to the cases in which $\bm=(0,\ldots,0)$ and $\bm \neq (0,\ldots,0)$.  Noting that $b_{(0,\ldots,0)} = 1$, we obtain 
$$\sum_{\bn \in \mathbb{N}_0^d} \Big(\E\prod_{j=1}^d  \frac{X_j^{n_j}}{n_j!}\Big) t_1^{n_1} \cdots t_d^{n_d} = \sum_{\bn \in \mathbb{N}_0^d} \bigg[\Big(\E\prod_{j=1}^p \frac{X_j^{n_j}}{n_j!}\Big) \cdot \Big(\E\prod_{j=p+1}^d \frac{X_j^{n_j}}{n_j!}\Big) + \delta_\bn\bigg] \prod_{j=1}^d t_j^{n_j},
$$
with $\delta_\bn \ge 0$.  Comparing the coefficients of the monomials $t_1^{n_1}\cdots t_d^{n_d}$, we obtain 
$$
\E\prod_{j=1}^d \Big(\frac{X_j^{n_j}}{n_j!} \Big) = \Big(\E\prod_{j=1}^p \frac{X_j^{n_j}}{n_j!}\Big) \cdot \Big(\E\prod_{j=p+1}^d \frac{X_j^{n_j}}{n_j!}\Big) + \delta_\bn \ge \Big(\E\prod_{j=1}^p \frac{X_j^{n_j}}{n_j!}\Big) \cdot \Big(\E\prod_{j=p+1}^d \frac{X_j^{n_j}}{n_j!}\Big),
$$
which yields \eqref{strong_gpi_gamma}, the strong form of the GPI. 
$\qed$

\bigskip

The following result was obtained by Russell and Sun \cite{Russell} by different methods.  In the context of Theorem \ref{thm_nonneg_corr}, the corollary follows from the well-known result that if $(X_1,\ldots,X_d) \sim N_d(0,\Sigma)$ then $\tfrac{1}{2} (X_1^2,\ldots,X_d^2) \sim {\hbox{Gamma}\hskip1pt}_d(\frac{1}{2},\Sigma)$.

\begin{corollary} \, {\rm{(Russell and Sun \cite{Russell})}}
	\label{cor_normal_strong_GPI}
	Let $(X_1,\ldots,X_d) \sim N_d(0,\Sigma)$, and suppose that there exists a sign matrix $S$ such that all off-diagonal elements of the matrix $S \Sigma S$ are nonnegative.  Then the strong form of the GPI \eqref{strong_gpi} holds for all even integers $n_1,\ldots,n_d$.  
\end{corollary}

\begin{remark}
	\label{remark_Isserlis_Wick}
	{\rm 
		An alternative approach to establishing Corollary \ref{cor_normal_strong_GPI} is by means of the classical Isserlis-Wick formula \cite{Kan}.  To see this, we write 
		\begin{align} 
		\E \prod_{j=1}^p X_j^{2m_j} &= \E\Big(\underbrace{X_1 \cdots X_1}_{2m_1 \, \rm{terms}} \cdot \underbrace{X_2 \cdots X_2}_{2m_2 \, \rm{terms}} \cdots  \underbrace{X_p \cdots X_p}_{2m_p \, \rm{terms}}\Big), \label{eq:reform2} \\
		\E \prod_{j=p+1}^d X_j^{2m_j} &= \E\Big(\underbrace{X_{p+1} \cdots X_{p+1}}_{2m_{p+1} \, \rm{terms}} \cdot \underbrace{X_{p+2} \cdots X_{p+2}}_{2m_{p+2} \, \rm{terms}} \cdots  \underbrace{X_d \cdots X_d}_{2m_d \, \rm{terms}}\Big), \label{eq:reform3} \\
		\intertext{and} 
		\E \prod_{j=1}^d X_j^{2m_j} &= \E\Big(\underbrace{X_1 \cdots X_1}_{2m_1 \, \rm{terms}} \cdot \underbrace{X_2 \cdots X_2}_{2m_2 \, \rm{terms}} \cdots  \underbrace{X_d \cdots X_d}_{2m_d \, \rm{terms}}\Big).\label{eq:reform1}
		\end{align}
		By the Isserlis-Wick formula, the expectations \eqref{eq:reform2}, \eqref{eq:reform3}, and \eqref{eq:reform1} can be written as a sum of products of the elements of $\Sigma_{11}$, $\Sigma_{22}$, and $\Sigma$, respectively.  Moreover, a simple inspection of the terms arising in the evaluation of \eqref{eq:reform2} and \eqref{eq:reform3} show that the product of those two expectations yields a collection of terms that are a subset of the terms arising from evaluation of \eqref{eq:reform1}.  Since we assume that $\sigma_{ij} \ge 0$ for all $i$ and $j$ then we obtain the strong form of the GPI.  
}\end{remark}

\begin{remark}
	\label{remark_thm_2_1_extended}
	{\rm We note that Theorem \ref{thm_nonneg_corr} can be extended further to distributions more general than the multivariate gamma distributions.  
	Consider mutually independent random vectors $Y_1,\ldots,Y_p \in \R^d$ such that, for all $i=1,\ldots,p$, $Y_i \sim {\hbox{Gamma}\hskip1pt}_d(\alpha_i,\Sigma_i)$ where all entries of the matrix $\Sigma_i$ are nonnegative.  Denote by $Y_{ij}$ the $j$-th component of $Y_i$, and define a random vector $(Z_1,\ldots,Z_d)$ by $Z_j = \sum_{i=1}^p Y_{ij}$, $j=1,\ldots,d$.  Then it is straightforward to show that the moment-generating function of $(Z_1,\ldots,Z_d)$ is
		$$
		\prod_{i=1}^p \det(I_d-\Sigma_i T)^{-\alpha_i} = \exp\Big(\sum_{i=1}^p \alpha_i \sum_{n=1}^\infty \frac{1}{n} \tr[(\Sigma_i T)^n]\Big).
		$$
	 Exploiting the additivity of the traces and using similar arguments as in the proof of Theorem \ref{thm_nonneg_corr} yields an analogous theorem for the vector $(Z_1,\ldots,Z_d)$.}
\end{remark}

\section{Positive upper orthant dependence and the GPI}
\label{sec_gpi}

In this section, we investigate the validity of the inequalities \eqref{weak_gpi} and \eqref{strong_gpi} without making specific assumptions on the distribution of the marginals of $\bX$. As already pointed out in the introduction, \eqref{strong_gpi} is valid if $|\bX|$, the vector of absolute values, is associated. It is also clear that \eqref{strong_gpi} holds when $|\bX|$ is {\it weakly associated} \cite{zheng}.

In this section, we show that \eqref{strong_gpi} follows from the notion of strong positive upper orthant dependence (SPUOD), which has been shown to be strictly weaker than weak association \cite{zheng}. Moreover the weak form \eqref{weak_gpi} follows from the notion of positive upper orthant dependence (PUOD).

Let us recall \cite{Dharma} that a random vector $V = (V_1,\ldots,V_d) \in \R^d$ is said to be {\it positively upper orthant dependent} (PUOD) if 
$$
\Pr(V_1 \ge t_1,\ldots,V_d \ge t_d) \ge \prod_{j=1}^d \Pr(V_j \ge t_j) 
$$
for all $t_1,\ldots,t_d \in \R$.  We will also say that the vector $V$ is \textit{strongly positively upper orthant dependent} (SPUOD) if 
$$
\Pr(V_1 \ge t_1,\ldots,V_d \ge t_d) \ge \prod_{j \in I} \Pr(V_j \ge t_j) \cdot \prod_{j \in I^c} \Pr(V_j \ge t_j)
$$
for all $I \subset 1,\ldots,d$.  

We begin with a result which, in the one-dimensional case, is classical in the literature on the statistical analysis of survival data.  

\begin{lemma}
	\label{lemma_survival_mean}
	Let $Y = (Y_1,\ldots,Y_d)$ be a random vector with nonnegative components $Y_1,\ldots,Y_d$ and such that $\E(Y_1 \cdots Y_d) < \infty$.  Then 
	\begin{equation}
	\label{survival_mean}
	\E(Y_1 \cdots Y_d) = \int_0^\infty \cdots \int_0^\infty \Pr(Y_1 \ge t_1,\ldots,Y_d \ge t_d) \dd t_1 \cdots \dd t_d.
	\end{equation}
\end{lemma}

\Pro
For completeness, we provide a direct proof; cf., \cite{Liu}.  For $y \ge 0$, let $\chi_y$ be the indicator function of the interval $[0,y]$; i.e., $\chi_y(t) = 1$ if $0 \le t \le y$, and $\chi_y(t) = 0$ if $t > y$.  Then 
$$
\int_0^\infty \chi_y(t) \dd t = y,
$$
and it follows by an application of Fubini's theorem that 
\begin{align}
\label{E_of_Ys}
\E(Y_1 \cdots Y_d) &= \E \prod_{j=1}^d \int_0^\infty \chi_{Y_j}(t_j) \dd t_j \nonumber \\
&= \int_0^\infty \cdots \int_0^\infty \E \prod_{j=1}^d \chi_{Y_j}(t_j) \dd t_j.
\end{align}
It is trivial that 
$$
\prod_{j=1}^d \chi_{Y_j}(t_j) = 
\begin{cases}
1, & \hbox{ if } Y_1 \geq t_1,\ldots,Y_d \geq t_d \\
0, & \hbox{ otherwise}
\end{cases};
$$
therefore 
$$
\E \prod_{j=1}^d \chi_{Y_j}(t_j) = \Pr(Y_1 \geq t_1,\ldots,Y_d \geq t_d).
$$
Substituting the latter result into \eqref{E_of_Ys}, we obtain \eqref{survival_mean}.  
$\qed$

\medskip

\begin{theorem}
	\label{thm_gpi}
	Let $(X_1,\ldots,X_d)$ be a random vector such that $\E(|X_1|^{n_1} \cdots |X_d|^{n_d}) < \infty$ for fixed exponents $n_1,\ldots,n_d \ge 0$.  
	
	If $|\bX|$, the vector of absolute values of $(X_1,\ldots,X_d)$, is PUOD then there holds the weak form of the GPI, 
	\begin{equation} 
	\label{eq:puodweak}
	\E \Big(\prod_{j=1}^d |X_j|^{n_j}\Big) \ge \prod_{j=1}^d \E\big(|X_j|^{n_j}\big).
	\end{equation}
	
	If $|\bX|$ is SPUOD then the strong form of the GPI holds, i.e., for any $I \subset \{1,\ldots,d\}$, 
	\begin{equation} \label{eq:puodstrong}
	\E \Big(\prod_{j=1}^d |X_j|^{n_j}\Big) \ge \E \Big(\prod_{j \in I} |X_j|^{n_j}\Big) \cdot \E\Big(\prod_{j \in I^c} |X_j|^{n_j}\Big).
	\end{equation}
\end{theorem}

\Pro
Suppose that $|\bX|$ is PUOD.  Replacing each $t_j$ by $t_j^{1/n_j}$ and simplifying the various inequalities on the $X_j$, $j=1,\ldots,d$, we obtain  
$$
\Pr(|X_1|^{n_1} \ge t_1,\ldots,|X_d|^{n_d} \ge t_d) \geq \prod_{j=1}^d \Pr(|X_j|^{n_j} \ge t_j)
$$
for all $t_1,\ldots,t_d \ge 0$.  Integrating both sides of this inequality with respect to $t_1,\ldots,t_d$ over $\R^d_{\ge 0}$ then, by applying Lemma \ref{lemma_survival_mean}, we obtain \eqref{eq:puodweak}.

The strong form of the GPI \eqref{eq:puodstrong} can be derived analogously starting from the assumption that $|\bX|$ is SPUOD.  
$\qed$

\section{The strong form of the GPI for negative exponents}
\label{sec_negative_exponents}

The strong form of the GPI \eqref{strong_gpi} for the case in which all exponents are negative was proved by Wei \cite{Wei}. We now derive this result succinctly by an application of the Gaussian correlation inequality \cite{Royen} and the method of integrating the multivariate survival function, as applied earlier in Section \ref{sec_gpi}.

\begin{proposition}
	\label{prop_negative_exp}
	Suppose that $(X_1,\ldots,X_d) \sim N_d(0,\Sigma)$ and that $n_1,\ldots,n_d \in (0,1)$.  Then 
	\begin{equation}
	\label{s_gpi_negative_exp}
	\E \Big(\prod_{j=1}^d |X_j|^{-n_j}\Big) \ge \E \Big(\prod_{j \in I} |X_j|^{-n_j}\Big) \cdot \E \Big(\prod_{j \in I^c} |X_j|^{-n_j}\Big),
	\end{equation}
	for all $I \subset \{1,\ldots,n\}$
\end{proposition}

\Pro
Without loss of generality, we can assume that $I = \{1,\ldots,p\}$.  We note that the conditions on $n_1,\ldots,n_d$ are necessary to ensure that the moments in \eqref{s_gpi_negative_exp} are finite.  

For $t_1,\ldots,t_d > 0$, we apply the Gaussian correlation inequality \cite{Royen} to obtain 
\begin{align*}
\Pr(| & X_1|^{-n_1} \ge t_1,\ldots,|X_d|^{-n_d} \ge t_d) \\
&= \Pr(|X_1| \le t_1^{-1/n_1},\ldots,|X_d| \le t_d^{-1/n_d}) \\
&\ge \Pr(|X_1| \le t_1^{-1/n_1},\ldots,|X_p| \le t_p^{-1/n_p}) \, \Pr(|X_{p+1}| \le t_{p+1}^{-1/n_{p+1}},\ldots,|X_d| \le t_d^{-1/n_d}) \\
&= \Pr(|X_1|^{-n_1} \ge t_1,\ldots,|X_p|^{-n_p} \ge t_p) \,  \Pr(|X_{p+1}|^{-n_{p+1}} \ge t_{p+1},\ldots,|X_d|^{-n_d} \ge t_d).
\end{align*}
Integrating the first and last terms of this inequality with respect to $t_1,\ldots,t_d \in (0,\infty)$, we obtain \eqref{s_gpi_negative_exp}.  
$\qed$

\bigskip

The argument used to prove Proposition \ref{prop_negative_exp} also establishes the novel finding that if any random vector $(Y_1,\ldots,Y_d)$ with (almost surely) positive components satisfies the Gaussian-type correlation inequality, 
$$
\Pr(Y_1 \le t_1,\ldots,Y_d \le t_d) \ge \Pr(Y_1 \le t_1,\ldots,Y_p \le t_p) \Pr(Y_{p+1} \le t_{p+1},\ldots,Y_d \le t_d)
$$
for all $t_1,\ldots,t_d > 0$, then 
$$
\E \Big(\prod_{j=1}^d Y_j^{-n_j}\Big) \ge \E \Big(\prod_{j=1}^p Y_j^{-n_j}\Big) \cdot \E \Big(\prod_{j=p+1}^d Y_j^{-n_j}\Big),
$$
for all $n_1,\ldots,n_d < 0$ such that the expectations exist, and for all $1 \le p \le d-1$.  In particular, the strong form of the GPI with negative exponents holds for the multivariate gamma distributions treated in \cite{Royen}.

\bigskip

\noindent
{\bf Acknowledgments}.  We are grateful to Fr\'ed\'eric Ouimet for drawing our attention to the article \cite{Ouimet} which motivated us to take another look at the GPI.

\vskip 0.25truein

\bibliographystyle{ims}

\end{document}